%% file: Header.tex
\documentclass[a4paper,abstracton]{scrartcl}
\usepackage{amssymb}
\usepackage{amsmath}
\usepackage{amsthm}
\usepackage{graphicx}
\usepackage{fancyhdr}              
\usepackage{lastpage}              
                                   
\usepackage[shortlabels]{enumitem}
\usepackage{esint}

\usepackage [english]{babel}

\newcommand{\setR}{\mathbb{R}}

\newcommand{\setZ}{\mathbb{Z}}
\newcommand{\setQ}{\mathbb{Q}}
\newcommand{\setN}{\mathbb{N}}

\newtheorem{definition}{Definition}
\newtheorem{theorem}[definition]{Theorem}
\newtheorem{lemma}[definition]{Lemma}
\newtheorem{corollary}[definition]{Corollary}

\newcommand{\abs}[1]{\left|#1\right|}

\newcommand{\floor}[1]{\left\lfloor #1\right\rfloor}

\newcommand{\set}[1]{\left\lbrace #1 \right\rbrace}
\newcommand{\card}[1]{\#\set{#1}}

\newcommand{\divides}{\mid}

\renewcommand{\pmod}[1]{(\text{mod } #1)}

\begin{document}
\title{Pairs of $k$-free Numbers, consecutive square-full Numbers}
\author{Thomas Reuss\\ 
		Mathematical Institute, University of Oxford\\
		reuss@maths.ox.ac.uk}
\date{\today}
\maketitle

\begin{abstract}
We consider the error term of the asymptotic formula for the number of pairs of $k$-free integers up to $x$. Our error term improves results by Heath-Brown, Brandes and Dietmann/Marmon. We then extend our results to $r$-tuples of $k$-free numbers and improve previous results by Tsang. Furthermore, we establish an error term for consecutive square-full integers. Finally, we will show that for all $\theta<3$ and for almost all $D$, the fundamental solution $\epsilon_D$ associated to the Pell equation $x^2-Dy^2=1$ satisfies $\epsilon_D> D^\theta$. This improves/recovers previous results by Fouvry and Jouve. The main tool of our work is the approximate determinant method. 
\end{abstract}

\renewcommand{\abstractname}{Acknowledgment}
\begin{abstract}
I am very grateful to my supervisor Roger Heath-Brown who introduced me to the determinant method and its applications. I thank him for many constructive and helpful comments on my work.\\
I would like to thank \'{E}tienne Fouvry and Florent Jouve for drawing my attention to their paper \cite{FouvryJouve} which resulted in Theorem \ref{thm:thmPell1} below.\\
I am also very grateful to the EPSRC\footnote{DTG reference number: EP/J500495/1} and to St. Anne's College, Oxford who are generously funding and supporting my doctorate's degree and this project.
\end{abstract}

\input{Introduction}
\input{DetMethod}
\input{Counting}
\input{Finishing}
\input{PairsOfKFree}
\input{FurtherImprovements}
\input{TouplesOfKFree}
\input{ConsecutiveSqFull}
\input{PellEquation}

\input{bibliography}

\end{document}

%% file: Introduction.tex
\section{Introduction}
\subsection*{Notation}
We write $x\sim X$ to say that $X<x<2X$ and we write $x\asymp X$ to say that there
exist positive constants A, B, independent of $X$, such that $AX \leq |x| \leq BX$.
\subsection*{The Main Theorem}
In this paper, we will present a generalization of the approximate determinant method developed by Heath-Brown in \cite{RHBsqfreen21}. More precisely, we will prove the following theorem:
\begin{theorem}\label{thm:thmfundamental}
Let $D,E,x>1$ and $\epsilon>0$. Let $k,l,h$ be integers such that
$1 \leq l<k$ and $h\neq 0$. Suppose $\mathcal{N}(x;D,E)$ is the number of elements in the set
\[
	\{(d,e,u,v)\in\setN^4: d\sim D,
	e\sim E, u\sim U, v\sim V, e^kv^l-d^ku^l=h\},
\]
where
\[
	U=\frac{x^{1/l}}{D^{k/l}},\quad\text{and}\quad V=\frac{x^{1/l}}{E^{k/l}}.
\]
Let $M>1$ be defined by
\[
	\log(M)=\frac{9}{8}\frac{\log(DE)\log(UV)}{\log x},
\]
and suppose that the following conditions are satisfied:
\begin{enumerate}
\item $\log(DE)\asymp\log(UV)\asymp\log(x)$.
\item $l\geq 2$, or $DE\gg_{k,l,h} x^{1/k}$.
\end{enumerate}
Then, if $x$ is large enough in terms of $\epsilon$,
\[
	\mathcal{N}(x;D,E)\ll_{\epsilon,k,l,h}
	x^\epsilon\min\{(DEM)^{1/2}+D+E,(UVM)^{1/2}+U+V\}.
\]
\end{theorem}
Theorem \ref{thm:thmfundamental} gives an upper bound for the number of integer points on the algebraic variety defined by
\begin{equation}\label{eq:eq1}
	e^kv^l-d^ku^l=h,
\end{equation}
where each of the variables $d,e,u,v$ are restricted to certain sizes $D,E,U,V$.
It is possible to generalize our result to algebraic varieties described by
$f_1(d,e)v^l+f_2(d,e)u^l=h$ where $f_1$ and $f_2$ are homogeneous polynomials of weighted degree. But for our applications, Theorem \ref{thm:thmfundamental} suffices.
The main strategy of the proof of Theorem \ref{thm:thmfundamental} will be to reduce the problem of counting integer points on the algebraic variety $e^kv^l-d^ku^l=h$ to counting rational points ``close'' to the curve $t=s^{k/l}$ where $t=v/u$ and $s=d/e$. The fundamental tool we use to tackle this counting problem will be the Determinant Method for which the interested reader should consult Heath-Brown \cite{cime}. 
Indeed, our proof of Theorem \ref{thm:thm1} is very similar in many stages to Heath-Brown \cite{RHBsqfreen21} where he derives an asymptotic formula for square-free values of the form $n^2+1$ up to $x$, say. In the following paragraphs, we will illustrate our applications of Theorem \ref{thm:thmfundamental}.

\subsection*{Pairs of $k$-free integers}
For $k\in\setZ_{\geq 2}$, we say that an integer $n$ is $k$-free if there is no prime $p$ such that $p^k\divides n$. By convention, a $2$-free integer is called square-free. It is an elementary fact that
\[
	S_k(x):=\card{n\leq x: n\text{ is } k\text{-free}}=\frac{x}{\zeta(k)}+O(x^{1/k}).
\]
We shall now consider the more general problem of deriving an asymptotic formula for pairs of $k$-free integers up to $x$, say.
More precisely, for integers $k\geq 2$ and $h\neq 0$, let $N_{k,h}(x)$ be the number of integers $n\leq x$ such that both $n$ and $n+h$ are $k$-free. Then, we have the following theorem:
\begin{theorem}\label{thm:thm1}
For all $\epsilon>0$ and all sufficiently large $x$, we have that
\[	
	N_{k,h}(x)=c_{k,h} x + O_{\epsilon,k,h}(x^{\omega(k)+\epsilon}),
\]
where
\[
	c_{k,h}=\prod_{p}\left(1-\frac{\rho_{k,h}(p)}{p^k}\right),
\]
with
\[
	\rho_{k,h}(p)= \begin{cases}
		2 & \text{if } p^k\nmid h\\
		1 & \text{if } p^k\divides h
	\end{cases},
\]
and
\[
	\omega(k)= \begin{cases}
	\frac{26+\sqrt{433}}{81}\approx 0.578 & \text{if } k=2\\
	\frac{169}{144k} & \text{for } k\geq 3.
	\end{cases}
\]
\end{theorem}
The problem of estimating $N_{k,h}(x)$ as in Theorem \ref{thm:thm1} is most interesting when $k$ is fairly small because $S_k(x)/x\rightarrow 1$ as $k\rightarrow\infty$. The trivial error term for Theorem \ref{thm:thm1} is $\omega(k)=2/(k+1)$ which can be obtained with elementary arguments (see for example Carlitz \cite{Carlitz}).\smallskip\\
By considering the Dirichlet series of $\zeta(s)/\zeta(ks)$, it is plausible to suspect that the error term of $S_k(x)$ cannot be improved by a proper power of $x$ below $O(x^{1/k})$ without assuming some quasi-Riemann-Hypothesis. By assuming uniformity in the error term of Theorem \ref{thm:thm1} with respect to $h\ll x$, one can show that an error term $O(x^\theta)$ in Theorem \ref{thm:thm1} implies an error term $O(x^\theta)$ in the formula for $S_k(x)$ and it is thus reasonable to assume that the error term of Theorem \ref{thm:thm1} cannot be improved below $O(x^{1/k})$ by a power of $x$ without radically new ideas.\smallskip\\
Heath-Brown \cite{RHBsquaresieve} considered the problem of Theorem \ref{thm:thm1} in the case $k=2$ and $h=1$. He obtained
$\omega(2)=7/11\approx 0.636$ improving the trivial value $\omega(2)=2/3$. Heath-Brown's approach uses exponential sums and the Square Sieve. It should be noted that Heath-Brown's method is actually uniform in $h\ll x$. Brandes \cite{Brandes} has generalized Heath-Brown's method for general $k$. Her argument yields the value $\omega(k)=14/(7k+8)$, which is the currently best available result for small $k$. Our error exponent $\omega(k)$ stated in Theorem \ref{thm:thm1} is such that $\omega(k)\leq 14/(7k+8)$ for all values of $k$. Dietmann and Marmon \cite{DietmannMarmon} have worked independently on improving the exponent $\omega(k)$ at the same time this paper was produced. Their method is similar to ours and they obtain the exponent $\omega(k)=14/(9k)$, which improves upon Brandes' bound for $k\geq 5$. Our error exponent is smaller than theirs for all $k\geq 2$. The proof of Theorem \ref{thm:thm1} will be an application of Theorem \ref{thm:thmfundamental} with $l=1$.

\subsection*{An asymptotic formula for $r$-tuples of $k$-free integers}
The second application of Theorem \ref{thm:thmfundamental} that we will present is the following:
\begin{theorem}\label{thm:thm2}
Let $k\geq 2$, $r\geq 2$ and $l_i(x)=a_ix+b_i\in\setZ[x]$ for $i=1,\ldots,r$ such that
$a_ib_j-a_jb_i\neq 0$ and $a_i\neq 0$ for all $i,j$ with $1\leq i,j\leq r$ and $i\neq j$. Then define
\[
	\rho(p)=\#\{n\pmod {p^k}: p^k\divides l_i(n)\text{ for some } i\},
\]
and let
\[
	c=\prod_{p}\left(1-\frac{\rho(p)}{p^k}\right).
\]
If $N(x)$ is the number of integers $n\leq x$ such that $l_1(n),\ldots,l_r(n)$ are all $k$-free. Then for any $\epsilon>0$ and any sufficiently large $x$ we have that 
\[
	N(x)=cx+O_\epsilon(x^{3/(2k+1)+\epsilon}).
\]
\end{theorem}
It should be noted that the implied constant in Theorem \ref{thm:thm2} depends on the choice of the $l_i$ and that the best error term in the sense of Theorem \ref{thm:thm2} available in the literature for $k=2$ was $O(x^{7/11+\epsilon})$ (See Tsang \cite{Tsang}). Tsang's proof uses a form of the Rosser-Iwaniec sieve and the version of Theorem \ref{thm:thm1} due to Heath-Brown. It should be noted that even though Tsang's error term is weaker than ours, his implied constants are uniform in $r$ and $\max_i ||l_i||$.

\subsection*{Consecutive square-full numbers}
We will prove the following theorem about consecutive square-full numbers. Recall that an integer $n$ is called square-full if, for all primes $p\divides n$, we have $p^2\divides n$.
\begin{theorem}\label{thm:thm3}
Let $N(x)$ be the number of integers $n\leq x$ such that both $n$ and $n+1$ are square-full. Then we have for all $\epsilon>0$ and sufficiently large $x$ that
\[
	N(x)\ll_\epsilon x^{29/100+\epsilon}.
\]
\end{theorem}
It can be shown that there are indeed infinitely many consecutive square-full numbers, for if $n$ and $n+1$ are square-full then so are $4n(n+1)$ and $4n(n+1)+1$. However, it follows from a simple application of the abc-conjecture that there are at most finitely many $n$ such that $n,n+1$ and $n+2$ are all square-full. The proof of Theorem \ref{thm:thm3} will be an application of Theorem \ref{thm:thmfundamental} with the variety given by $e^3v^2-d^3u^2=1$.

\subsection*{Size of the fundamental solution of Pell Equations}
Let $D>0$ be an integer which is not a square. Consider the Pell equation
\begin{equation}\label{eq:eqPell}
	T^2-DU^2=1
\end{equation}
as an equation in $(T,U)$. Write a solution of this equation as $\eta_D=T+U\sqrt{D}$.  For each fixed $D$, the Pell equation \eqref{eq:eqPell} has a unique solution called the fundamental solution $\epsilon_D$ which satisfies $\epsilon_D>1$ and 
\[
	\{\eta_D:\eta_D\text{ is a solution of \eqref{eq:eqPell}}\}
	=\{\pm\epsilon_D^n: n\in\setZ\}.
\]
It is strongly believed that $\epsilon_D$ is almost always much larger than $D$. In particular, one may study the quantity
\[
	\mathcal{S}(X,\alpha)=\#\{(D,\eta_D): D\sim X, 
	D\text{ is not a square}, 1<\eta_D\leq D^{1/2+\alpha}\}
\]
for $\alpha>1/2$. One can show that
\begin{equation}\label{eq:eqtrivFouvry}
	\mathcal{S}(X,\alpha)\ll X^{\min(\alpha+\epsilon,1)}.
\end{equation}
Fouvry and Jouve \cite{FouvryJouve} have shown that this estimate can be improved to
\[
	\mathcal{S}(X,\alpha)\ll_\epsilon X^{\frac{\alpha}{3}+\frac{7}{12}+\epsilon}.
\]
Their bound improves \eqref{eq:eqtrivFouvry} in the range $7/8< \alpha< 5/4$. We will
improve their result further. More precisely, we will prove the following result:
\begin{theorem}\label{thm:thmPell1}
\[
	\mathcal{S}(X,\alpha)\ll_\epsilon X^{\max\left(
	\frac{9}{16}\frac{\alpha}{\alpha+1/2}+\frac{1}{2}\min(1,\alpha),\frac{1}{2}+\frac{\alpha}{5}
	\right)}.
\]
\end{theorem}
Our result is stronger than \eqref{eq:eqtrivFouvry} in the range $5/8< \alpha< 5/2$.
We obtain the following corollary:
\begin{corollary}\label{cor:corPell}
	For all $0\leq\theta<3$, we have
	\[
		\#\{D: 1<D\leq X, D\text{ is not a square}, \epsilon_D\leq D^\theta\}=o_\theta(X).
	\]
\end{corollary}
This improves upon the result by Fouvry and Jouve \cite{FouvryJouve}. They have shown that one can take $\theta<7/4$ in Corollary \ref{cor:corPell}. Furthermore, Fouvry and Jove have proved a slightly weaker version of Corollary \ref{cor:corPell} in \cite{FouvryJouve2}. They showed that for all $\theta<3$, there is a positive proportion of fundamental discriminants $D$ such that $\epsilon_D >D^\theta$.
 For a more thorough discussion of the history and motivation behind Theorem \ref{thm:thmPell1} and Corollary \ref{cor:corPell}, the interested reader should consult their papers \cite{FouvryJouve}
 and \cite{FouvryJouve2}.

%% file: DetMethod.tex
\section{The Proof of Theorem \ref{thm:thmfundamental}}
We will now start to prove Theorem \ref{thm:thmfundamental}. From now on all implied constants may depend on $k,l$ and $h$. We may assume without loss of generality that $(ev,du)=1$, because if $p\divides (ev,du)$ then $p\divides h$. As a first step, we will set $t=v/u$ and $s=d/e$ so that the equation \eqref{eq:eq1} becomes
\[
	t^l=s^k+O\left(\frac{1}{E^kU^l}\right).
\]
Note that
\[
	t^l-s^k=(t-s^{k/l})\sum_{a+b=l-1}t^a(s^{k/l})^b.
\]
Using that $s\asymp D/E$ and $t\asymp V/U=(D/E)^{k/l}$, we may deduce that
\[
	\sum_{a+b=l-1}t^a(s^{k/l})^b\asymp \left(\frac{V}{U}\right)^{l-1}.
\]
Note that $E^k=x/V^l$ and hence
\begin{equation}\label{eq:ts2}
	t=s^{k/l}+O\left(\frac{1}{x}\frac{V}{U}\right).
\end{equation}
Thus, we have transformed our problem of counting integer points on a three-fold into a problem where we count rational points close to the curve $t=s^{k/l}$, where the sizes of the numerators and denominators of $s$ and $t$ are determined by $D$ and $E$. The determinant method seems to be stronger in counting-problems involving varieties of lower dimension so that dealing with a curve rather than a three-fold will provide the key saving in our proof.
In section \ref{sec:DetMethod} we will show how the determinant method allows us to subdivide the range of $s$ into intervals $I$ of equal length so that our problem transforms into counting rational points close to the curve $t=s^{k/l}$ where $s$
belongs to some particular interval $I$. In section \ref{sec:Counting}, we will then calculate the contribution of one such interval $I$ to $\mathcal{N}(x;D,E)$ and in section \ref{sec:Finishing} we will then add up all the contributions to finish the proof of the estimate in Theorem \ref{thm:thmfundamental}.
 
\subsection{Determinant Method}\label{sec:DetMethod}
Note that $s$ is of exact order $D/E$. Our plan is to pick an integer $M$ such that $\log x\ll \log M\leq\log x$ and split the range of $s$ into $O(M)$ intervals $I=(s_0,s_0(1+M^{-1})]$. For the rest of this section we will fix one such interval $I$
and consider solutions $(s,t)$ of \eqref{eq:ts2} with $s\in I$.
We label these solutions as $(s_1,t_1),\ldots,(s_J,t_J)$, say. Consider one such $(s_j,t_j)$. We can write $s_j=s_0(1+\alpha_j)$ where $0<\alpha_j\leq 1/M$: Note that $s_0^{k/l}\asymp V/U$ and hence
\[
	t_j=s_j^{k/l}+O\left(\frac{1}{x}\frac{V}{U}\right)
	=s_0^{k/l}\left((1+\alpha_j)^{k/l}+O\left(\frac{1}{x}\right)\right).
\]
Hence, after applying Taylor's Theorem with a suitable degree in $\alpha_j$ to $(1+\alpha_j)^{k/l}$, we can write
\begin{alignat*}{2}
	s_j &=s_0(1+\alpha_j) &\text{with } &0<\alpha_j\leq \frac{1}{M},\\
	t_j &=s_0^{k/l}(1+p(\alpha_j)+\beta_j)\quad &\text{with }
	&\beta_j\ll \frac{1}{x},
\end{alignat*}
where $p(\alpha_j)$ is a polynomial in $\alpha_j$ with no constant coefficient and with coefficients of size $O(1)$.
The next step is to choose positive integers $A$ and $B$ and to label the monomials $s^at^b$ with $a\leq A$ and $b\leq B$ as
$m_1(s,t),\ldots,m_H(s,t)$, where $H=(A+1)(B+1)$. Then one considers the $J\times H$ matrix $\mathcal{M}$
whose $(j,h)$-th entry is $m_h(s_j,t_j)$. We will show that the rank of $\mathcal{M}$ is strictly less than $H$ 
provided we choose $A,B$ and $M$ appropriately. This will enable us to deduce that there is a non-zero vector
$\mathbf{c}$ such that $\mathcal{M}\mathbf{c}=0$. If we now consider the polynomial $C_I(s,t)=C(s,t)=\sum_{h=1}^H c_h m_h(s,t)$ then we can see that $C(s_j,t_j)=0$
for all our solutions with $s_j\in I$. The vector $\mathbf{c}$ can be constructed from subdeterminants
of $\mathcal{M}$ which shows that $\mathbf{c}\in\setQ^H$ has rational entries. Note that the values $s_j$ and $t_j$ have numerators and denominators of size $\ll$ some power of $x$. By clearing out the common denominator of the coefficients of $C$ we may therefore assume
that $C$ has integer coefficients of size $\ll x^{\kappa(A,B)}$, say.\bigskip\\
To show that $\mathcal{M}$ has rank strictly less than $H$ we can assume that $H\leq J$ since otherwise this gets trivial.
Thus, it suffices to show that every $H\times H$ subdeterminant of $\mathcal{M}$ vanishes. Without loss of generality
it suffices to show that the determinant $\Delta$ coming from the first $H$ rows and columns of $\mathcal{M}$ vanishes.
The $j$-th row of $\mathcal{M}$ has entries with common denominator $e_j^Au_j^B$ which implies that
\begin{equation}\label{eq:eqProdDelta}
	\left(\prod_{j\leq H} e_j^A u_j^B\right)\Delta\in\setZ.
\end{equation}
Observe that $\prod_{j\leq H} e_j^A u_j^B\ll E^{AH}U^{BH}$ and hence if we can show
\[
	|\Delta|\ll_{A,B} E^{-AH}U^{-BH},
\]
with a suitable implied constant, then the integer in \eqref{eq:eqProdDelta} has to be an integer strictly less than 1 which implies $\Delta=0$. We substitute $s_j=s_0(1+\alpha_j)$ and $t_j =s_0^{k/l}(1+p(\alpha_j)+\beta_j)$ so that $\mathcal{M}$ has entries 
\[
	s_0^{a+bk/l}(1+\alpha_j)^a(1+p(\alpha_j)+\beta_j)^b.
\]
Hence
\[
	\Delta=\left(\prod_{a=0}^A\prod_{b=0}^B s_0^{a+bk/l}\right)\Delta_1=s_0^{\frac{1}{2}H(A+Bk/l)}\Delta_1,
\]
where $\Delta_1$ is the determinant of the generalized Vandermonde matrix with its entries
being polynomials in $\alpha_j$ and $\beta_j$ and coefficients of size $O_{A,B}(1)$. Note that we have
\[
	|\alpha_j|\leq X_1^{-1}\text{ and } |\beta_j|\leq X_2^{-1},
\]
where $X_1$ is of exact order $M$ and
$X_2$ is of exact order $x$. In particular, note that $\log X_1$ and $\log X_2$ are both of exact order $\log x$.
We now order the monomials $X_1^{-a}X_2^{-b}$ decreasing in size, $1=M_0,M_1,\ldots,M_H,\ldots$ say. Then by Lemma 3 in Heath-Brown
\cite{RHBdiffthreekpow}, we may bound $\Delta_1$ and hence $\Delta$ as follows:
\[
	\Delta=s_0^{\frac{1}{2}H(A+Bk/l)}\Delta_1 \ll_{A,B} (D/E)^{\frac{1}{2}H(A+Bk/l)}\prod_{h=1}^H M_h.
\]
Let $M_H=W^{-1}$. Then $X_1^{-a}X_2^{-b}\geq M_H$ if and only if
\begin{equation}\label{eq:ineqforW}
	a\log X_1+b\log X_2\leq\log W.
\end{equation}
The number of pairs $(a,b)$ which satisfy this inequality is
\[
	\frac{1}{2}\frac{(\log W)^2}{(\log X_1)(\log X_2)}
	+O\left(\frac{\log W}{\log x}\right)+O(1).
\]
This number must be equal to $H$ which gives
\[
	|2H(\log X_2)(\log X_1)-(\log W)^2|\ll (\log W)(\log x)+(\log x)^2,
\]
and hence
\[
	|\sqrt{2H(\log X_2)(\log X_1)}-\log W|\ll (\log x)\frac{(\log W)+(\log x)}{\sqrt{2H(\log X_2)(\log X_1)}+\log W}\ll\log x,
\] 
since $(\log X_2)(\log X_1)\asymp (\log x)^2$. Therefore
\begin{equation}\label{eq:eqforlogW}
	\log W = \sqrt{2H(\log X_2)(\log X_1)}+O(\log x).
\end{equation}
Next, observe that
\begin{align*}
	\log\prod_{h=1}^H M_h &=-\sum_{a,b}(a\log X_1+b\log X_2)\\
	&= -\frac{1}{3}\frac{(\log W)^3}{(\log X_1)(\log X_2)}
	  +O\left(\frac{(\log W)^2}{\log x}\right)+O(\log x),
\end{align*}
where the summation is subject to the inequality \eqref{eq:ineqforW}. 
Using the value for $\log W$ in \eqref{eq:eqforlogW}, we get
\[
	\log\prod_{h=1}^H M_h = -\frac{2\sqrt 2}{3}H^\frac{3}{2}\sqrt{(\log X_1)(\log X_2)}+O(H\log x).
\]
This shows that
\[
	\log|\Delta|\leq O_{A,B}(1)+\frac{1}{2}H(A+B\frac{k}{l})\log\frac{D}{E} -\frac{2\sqrt 2}{3}H^\frac{3}{2}\sqrt{(\log X_1)(\log X_2)}+O(H\log x).
\]
Thus, to get our required bound on $\log|\Delta|$ it suffices to show that
\[
	A\log E+B\log U+\frac{1}{2}(A+B\frac{k}{l})\log\frac{D}{E}+C_1(A,B)+C_2\log x\leq\frac{2\sqrt 2}{3}H^\frac{1}{2}\sqrt{(\log X_1)(\log X_2)},
\]
where $C_1(A,B)$ is some constant depending on $A$ and $B$ and $C_2$ is an absolute constant. Note that $AB\leq H$ and
\[
	A\log E+B\log U+\frac{1}{2}(A+B\frac{k}{l})\log\frac{D}{E}=\frac{A}{2}\log(DE)+\frac{B}{2}\log(UV).
\]
Hence, it suffices to show
\[
	\frac{1}{2}A\log(DE)+\frac{1}{2}B\log(UV)+C_1(A,B)+C_2\log x\leq\frac{2\sqrt 2}{3}(AB)^\frac{1}{2}\sqrt{(\log X_1)(\log X_2)}.
\]
We will optimize the error by choosing $A=\floor{\frac{B\log(UV)}{\log (DE)}}$.
The first assumption of Theorem \ref{thm:thmfundamental} implies that $\nu:=\frac{\log(UV)}{\log (DE)}\asymp 1$. In particular, if $B\geq 2/\nu$
then we may apply the inequalities $\frac{1}{2}B\nu\leq\floor {B\nu}\leq B\nu$ for $B\geq 2/\nu$ to get that $A\asymp B$. Analyzing Lemma 3 of
Heath-Brown \cite{RHBdiffthreekpow} which produces our constant $C_1(A,B)$, we can see that $C_1(A,B)$ can be replaced by some $C_3(B)$, say. Now observe that $\frac{1}{2}A\log(DE)+\frac{1}{2}B\log(UV)\leq B\log(UV)$ and that
\begin{align*}
	     \frac{2\sqrt 2}{3} (AB)^\frac{1}{2}\sqrt{(\log X_1)(\log X_2)}
	\geq & \frac{2\sqrt 2}{3} B\nu^\frac{1}{2}\left(1-\frac{1}{B\nu}\right)^\frac{1}{2}\sqrt{(\log X_1)(\log X_2)} \\
	\geq & \frac{2\sqrt 2}{3} B\nu^\frac{1}{2}\left(1-\frac{1}{B\nu}\right)\sqrt{(\log X_1)(\log X_2)} \\
	=    & \frac{2\sqrt 2}{3} B\sqrt{\frac{\log(UV)}{\log(DE)}}\sqrt{(\log X_1)(\log X_2)}+O(\log x)
\end{align*}
since $(\log X_1)(\log X_2)\asymp (\log x)^2$. It therefore suffices if we have
\[
	B\log(UV)+C_3(B)+C_2\log x\leq \frac{2\sqrt 2}{3} B\sqrt{\frac{\log(UV)}{\log(DE)}}\sqrt{(\log X_1)(\log X_2)}.
\]
Let $\delta>0$ be arbitrary and pick $X_1$ such that
\[
	\frac{2\sqrt 2}{3} \sqrt{\frac{\log(UV)}{\log(DE)}}\sqrt{(\log X_1)(\log X_2)}\geq (1+\delta)\log(UV),
\]
whence it suffices if we have $B\log(UV)+C_3(B)+C_2\log x\leq B(1+\delta)\log(UV)$ which holds if and only if
$B\delta\log(UV)\geq C_3(B)+C_2\log x$. Pick 
\[
	B=B(\delta)\geq\frac{2C_2\log x}{\delta\log(UV)}
\]
so that $C_2\log x\leq \frac{1}{2}\delta B\log(UV)$ and pick $x$ large enough in terms of $\delta$ so that
$C_3(B)\leq\frac{1}{2}\delta B\log(UV)$. Thus, we have shown that if we pick $B$ (and hence $A$) and $x$ large enough
in terms of $\delta$ and if we can pick a suitable $X_1$ as above then $\Delta=0$. Observing that $X_1\gg M$ and
$X_2\gg x$, we may rewrite the condition with a redefined $\delta$ as
\[
	\log M\geq \frac{9}{8}(1+\delta)\frac{\log(UV)\log(DE)}{\log x}.
\]
In summary, we get the following lemma:
\begin{lemma}\label{lemma1} Let $\delta>0$, and assume that an integer $M$ satisfies $\log x\ll \log M\leq \log x$ and
	\[
		\log M\geq \frac{9}{8}(1+\delta)\frac{\log(DE)\log(UV)}{\log x}.
	\]
	Then for any interval $I=\left(s_0,s_0(1+M^{-1})\right]$, there exists a non-zero integer polynomial
	$C_I(s,t)$ of total degree $O_\delta(1)$ and coefficients of size $O(x^\kappa)$ where $\kappa=\kappa(\delta)$
	and such that $C_I(s,t)=0$ for all solutions of
	\eqref{eq:ts2} with $s\in I$.
\end{lemma}

%% file: Counting.tex
\subsection{Counting Solutions}\label{sec:Counting}
Next, we want to estimate how much a fixed interval $I$ contributes to $\mathcal{N}(x;D,E)$. As in \cite{RHBsqfreen21}, we may assume that $C_I$ is absolutely irreducible. By clearing the denominators of $C_I(s,t)=0$, we may rewrite the equation in the form
\begin{equation}\label{eq:eqF}
	F(d,e;v,u)=0.
\end{equation}
In summary, we have created an auxiliary polynomial $F$ which is bi-homogeneous of degree $(a,b)$, say, has a total degree of size $O_\delta(1)$ and coefficients of size $O(x^{\kappa(\delta)})$ and it vanishes at the solutions $(s,t)$ of \eqref{eq:ts2} with $s\in I$.
The condition $s\in I$ gives $|d-es_0|\ll D/M$. We also have $|e|\ll E$. It is convinient to define the linear map
$T:\setR^2\rightarrow\setR^2$ by
\[
	T(x_1,x_2)=\left(\frac{M}{D}(x_1-x_2s_0),\frac{1}{E} x_2\right).
\]
This linear map defines a lattice 
\[
	\Lambda=\lbrace T(x_1,x_2):(x_1,x_2)\in\setZ^2\rbrace
\]
of determinant $\det(\Lambda)=M/(ED)$. Consider the rectangle
\[
	R=\{(\alpha_1,\alpha_2): |\alpha_1|\ll 1,|\alpha_2|\ll 1\},
\]
with suitable implied constants so that $s=d/e\in I$, implies $T(d,e)\in\Lambda\cap R$. Thus, we will now count points falling into $\Lambda\cap R$.\\ 
Let $\mathbf{g}^{(1)}$ be the shortest non-zero vector in $\Lambda$ and let $\mathbf{g}^{(2)}$ be the shortest vector in $\Lambda$
not parallel to $\mathbf{g}^{(1)}$. Then $\mathbf{g}^{(1)}, \mathbf{g}^{(2)}$ will be a basis for $\Lambda$. Moreover,
$\lambda_1\mathbf{g}^{(1)}+\lambda_2\mathbf{g}^{(2)}\in R$ implies $|\lambda_1 \mathbf{g}^{(1)}|\ll 1$ and $|\lambda_2 \mathbf{g}^{(2)}|\ll 1$.
By defining $L_i$ to be a suitable constant times
$|\mathbf{g}^{(i)}|^{-1}$ for $i=1,2$, we may write $|\lambda_i|\leq L_i$. Note that $|\mathbf{g}^{(1)}|\leq |\mathbf{g}^{(2)}|$ and 
$|\mathbf{g}^{(1)}||\mathbf{g}^{(2)}|\ll\det(\Lambda)=M/(DE)$ from which we may conclude $L_2\ll L_1$ and $L_1L_2\gg \frac{DE}{M}$.
Next set
\[
	\mathbf{h}^{(i)}=  (\frac{D}{M} \mathbf{g}_1^{(i)}+s_0 E\mathbf{g}_2^{(i)},E\mathbf{g}_2^{(i)}).
\]
Then $\mathbf{h}^{(1)}$ and $\mathbf{h}^{(2)}$ will form a basis for $\setZ^2$ and if
$\mathbf{x}=\lambda_1\mathbf{h}^{(1)}+\lambda_2\mathbf{h}^{(2)}$ is in the region given by $|x_1-x_2s_0|\ll D/M$ and
$|x_2|\ll E$ then $\lambda_1\mathbf{g}^{(1)}+\lambda_2\mathbf{g}^{(2)}\in R$. Thus, after a change of basis we may
replace $(d,e)$ by $(\lambda_1,\lambda_2)$ with $|\lambda_i|\leq L_i$ where $L_2\ll L_1$ and $L_1L_2\gg (DE)/M$.\bigskip\\
Our value $s_0$ is of the form $s_0=\frac{x_3}{M}\frac{D}{E}$ where $x_3$ is an integer of exact order $M$. 
Define $x_4:=\lfloor x_3^{k/l}M^{1-{k/l}}\rfloor$ and $t_0:=\frac{x_4}{M}\frac{V}{U}$ so that $s_0^{k/l}=t_0+O\left(\frac{1}{M}\frac{V}{U}\right)$.
Observe that
\[
	t=s_0^{k/l}+s_0^{k/l}p(\alpha)+s_0^{k/l}\beta=s_0^{k/l}+
	O\left(\frac{1}{M}\frac{V}{U}\right)+
	O\left(\frac{1}{x}\frac{V}{U}\right)
	=t_0+O\left(\frac{1}{M}\frac{V}{U}\right),
\]
since $M\leq x$.
This leads to the conditions $|v-u t_0|\ll \frac{V}{M}$ and $|u|\ll U$ and hence we can analogously to the
above argument replace $(v,u)$ by $(\tau_1,\tau_2)$ say where $|\tau_i|\leq T_i$ with $T_2\ll T_1$ and $T_1T_2\gg\frac{UV}{M}$.
These substitutions convert equation \eqref{eq:eq1} into an equation
\begin{equation}\label{eq:eqG0}
	G_0(\lambda_1,\lambda_2;\tau_1,\tau_2)= h,
\end{equation}
say where $G_0$ is bi-homogeneous of degree $(k,l)$. Similarly, equation \eqref{eq:eqF} will turn into an
equation of the form
\begin{equation}\label{eq:eqG1}
	G_1(\lambda_1,\lambda_2;\tau_1,\tau_2)=0,
\end{equation}
where $G_1$ is bi-homogeneous of degree $(a,b)$ and satisfies the same conditions as $F$. Also, from the above argument it is clear that 
the vectors $(\lambda_1,\lambda_2)$ and $(\tau_1,\tau_2)$ are primitive. For example, if $p\divides (\lambda_1,\lambda_2)$ then $p\divides (d_1,e_1)=1$.\bigskip\\
If $a=0$ then \eqref{eq:eqG1} determines $O_\delta(1)$ pairs $(\tau_1,\tau_2)$ each of which gives a pair $(u,v)$. The number of pairs $(d,e)$ corresponding to such a $(u,v)$ is 
\begin{equation}\label{eq:eqtrivialde}
	\#\{(d,e): d\asymp D, e\asymp E, e^k v^l-d^ku^l=h\}\ll x^\epsilon.
\end{equation}
This follows since $e^k v^l-d^ku^l=h$ is a Thue equation for $k\geq 3$ which has only $O_{k,l,h}(1)$ solutions (Thue \cite{Thue}). For $k=2$, the estimate \eqref{eq:eqtrivialde} follows from Estermann \cite{Estermann}. Thus, the case $a=0$ contributes $O_\delta(x^\delta)$ to
$\mathcal{N}(x;D,E)$.\bigskip\\
If $b=0$ then \eqref{eq:eqG1} determines $O_\delta(1)$ pairs $(\lambda_1,\lambda_2)$ each of which gives a pair $(d,e)$. The number of pairs $(u,v)$ corresponding to such a pair $(d,e)$ is
\begin{equation*}
	\#\{(u,v): u\asymp U, v\asymp V, e^k v^l-d^ku^l=h\}\ll x^\epsilon.
\end{equation*}
If $l=2$ then this follows again by Estermann's result and if $l\geq 3$ the estimate follows from the above result about Thue equations. If $l=1$ then the equation $e^kv-d^ku=h$ is a linear Diophantine equation in $(u,v)$ and thus, for each pair $(d,e)$, the number of solutions $(u,v)$ to \eqref{eq:eq1} is
\[
	\ll\frac{V}{D^k}+1=\frac{x}{(DE)^k}+1\ll 1,
\]
where the last estimate follows from the assumption $DE\gg x^{1/k}$ of the theorem. Hence, the case $b=0$ contributes $O_\delta(1)$ to $\mathcal{N}(x;D,E)$.\bigskip\\
If $a\geq 2$ then Lemma 2 of Heath-Brown \cite{RHBsqfreen21} gives us that \eqref{eq:eqG1} has $O_{\epsilon,\delta}(T_1^{1+\epsilon}x^{\kappa\epsilon})$
solutions. (recall that the coefficients of $G_1$ are of order $x^\kappa$). Picking $\epsilon$ small enough in terms of $\delta$,
we may assume that this is $O_{\delta}(T_1^{1+\delta}x^{\delta})$.
Each solution of \eqref{eq:eqG1} produces at most one solution of \eqref{eq:eq1}. So in the case $a\geq 2$ the contribution of
$I$ to $\mathcal{N}(x;D,E)$ is $O_{\delta}(T_1^{1+\delta}x^{\delta})$.
Similarly, if $b\geq 2$ then the contribution of $I$ to $\mathcal{N}(x;D,E)$ is $O_{\delta}(L_1^{1+\delta}x^{\delta})$.\bigskip\\
If $a=1$ then \eqref{eq:eqG1} can be written as
\[
	\lambda_1 G_{11}(\tau_1,\tau_2)+\lambda_2 G_{12}(\tau_1,\tau_2)=0,
\]
and hence $q\lambda_1=G_{12}(\tau_1,\tau_2)$ and $q\lambda_2=-G_{11}(\tau_1,\tau_2)$ for some integer $q$.
We define two polynomials $g_1,g_2\in\setZ[x]$ by $G_{1i}(\tau_1,\tau_2)=\tau_1^bg_i(\tau_2/\tau_1)$ for $i=1,2$.
Observe that $g_1$ and $g_2$ must be coprime since $G_1$ is absolutely irreducible. Hence, by the Euclid's Algorithm there
exist polynomials $h_1,h_2\in\setZ[x]$ and an integer $H$ such that
\[
	g_1h_1+g_2h_2=H
\]
where $H=O(x^\kappa)$. Evaluating this equation at $\tau_2/\tau_1$ we may deduce that $q\divides H\tau_1^{K}$ where $K$ is some integer. 
and similarly we may conclude that $q\divides \tilde H\tau_2^K$ where $\tilde H=O(x^\kappa)$.
But $\tau_1$ and $\tau_2$ are coprime. Thus, we may deduce that we have $O_\delta(x^\delta)$ choices for $q$.
Each value of $q$ gives us a value of $\lambda_1$ and $\lambda_2$ in terms of $\tau_1$ and $\tau_2$ which we may
substitute into \eqref{eq:eqG0} to get a Thue equation of the form $G_3(\tau_1,\tau_2)=h q^k$, say. 
This equation gives $O(T_1)$ possible pairs $\tau_1,\tau_2$ which shows that the case $a=1$ contributes
at most $O_\delta(x^\delta T_1)$ to $\mathcal{N}(x;D,E)$. Similarly, we may deduce that
the case $b=1$ contributes at most $O_\delta(x^\delta L_1)$ to $\mathcal{N}(x;D,E)$.
In summary, we obtain the following:
\begin{lemma}\label{lemma3} For any $\delta>0$, the contribution of the solutions $(s,t)$ with $s\in I$ to $\mathcal{N}(x;D,E)$ is
$O_\delta(x^\delta\min(L_1^{1+\delta},T_1^{1+\delta}))$.
\end{lemma}

%% file: Finishing.tex
\subsection{Completion of the proof of Theorem \ref{thm:thmfundamental}}\label{sec:Finishing}
In the previous section, we calculated the contribution of each interval $I$ to $\mathcal{N}(x;D,E)$. It remains to sum up the contribution of the various intervals. To proceed, we write $\mathbf{g}^{(1)}$ from the previous section as
\[
	\mathbf{g}^{(1)}=\left((M/D)(x_1-x_2s_0),(1/E)x_2\right).
\]
Recall that $L_1$ was defined to be a suitable multiple of 
$|\mathbf{g}^{(1)}|^{-1}$. This gives 
\[
	L_1(x_1-x_2s_0)\ll D/M,\quad\text{and}\quad L_1x_2\ll E.
\]
Recall that we produce the intervals $I=(s_0,s_0+\frac{1}{M}\frac{D}{E}]$ by taking $s_0=x_3\frac{1}{M}\frac{D}{E}$ for integers $x_3\asymp M$. Hence,
the number of intervals $I$ for which $L<L_1\leq 2L$ is at most the number of triples $(x_1,x_2,x_3)\in\setZ^3$
for which $\gcd(x_1,x_2)=1$ and
\[
	x_2x_3=\frac{ME}{D}x_1+O\left(\frac{E}{L}\right),\quad x_2\ll \frac{E}{L},\quad x_3\asymp M.
\]
Recalling that $L_1\gg L_2$ and $L_1L_2\gg (DE)/M$, we can deduce that $L\gg (DE/M)^{1/2}$.\bigskip\\
If $x_2=0$ then $x_1=\pm 1$ since $(x_1,x_2)=1$. Hence,
$L\asymp |\mathbf{g}^{(1)}|^{-1}=D/M$. Note that $x_3\neq 0$. Thus, there are $O(M)$ choices for $x_3$. Hence, the number of intervals $I$ for which $x_2=0$ and $L<L_1\leq 2L$ is $O(M)$. By Lemma \ref{lemma3}, each such interval contributes at most $O(x^\delta L^{1+\delta})$ to $\mathcal{N}(x;D,E)$ which gives a total contribution of $O(x^\delta D)$ corresponding to these intervals.\bigskip\\ 
Next consider the cases when $x_1=0$. In this case $x_2=\pm 1$ because $(x_1,x_2)=1$, and hence $L\asymp E/M$. And thus, a similar argument as above shows that the intervals $I$ for which $x_1=0$ and $L<L_1\leq 2L$ contribute a total of $O(x^\delta E)$ to $\mathcal{N}(x;D,E)$.\bigskip\\
We are therefore left with the case when $x_1$ and $x_2$ are both non-zero.
We can now see that $|\mathbf{g}^{(1)}|\gg E^{-1}$ which implies that
$L\ll E$. The above conditions on $(x_1,x_2,x_3)$ imply that $x_1\ll D/L$.
Thus, there are $O(D/L)$ choices for $x_1$ and each $x_1$ produces $O(E/L)$ choices for the product $x_2x_3$. And since $x_2x_3\neq 0$, a divisor function estimate shows that each value of $x_2x_3$ arises from at most $O_\delta(x^\delta)$ pairs $(x_2,x_3)$. Thus,
there are $O_\delta(x^\delta DE/L^2)$ intervals $I$ so that $L_1$ is of exact order $L$. Each interval contributes
$O_\delta(x^\delta L^{1+\delta})$ by Lemma \ref{lemma3}. Note that $L^\delta\ll E^\delta\ll x^\delta$ and hence we
get a total contribution of $O_\delta(x^{\delta}DE/L)$ corresponding to the intervals with $x_1x_2\neq 0$. By dyadic subdivision of the range of 
$L\gg (DE/M)^{1/2}$ we can conclude that
\[
	\mathcal{N}(x;D,E)\ll_\delta x^\delta ((DEM)^{1/2}+D+E).
\] 
Similarly, we may consider the $(x_1,x_2)$ lattice from the previous section corresponding to $(v,u)$ to get
$T_1(x_1-x_2t_0)\ll V/M$ and $T_1x_2\ll U$. Recall that 
$t_0=x_4\frac{1}{M}\frac{V}{U}$ where $x_4=\floor{\frac{x_3^k}{M^{k-1}}}$. We can see that $x_3\asymp M$ implies
$x_4\asymp M$ for large enough $x$. Hence by a completely analogous argument to the above we can deduce that
\[
	\mathcal{N}(x;D,E)\ll_\delta x^\delta ((UVM)^{1/2}+U+V).
\] 
Thus, we have established the bound for $\mathcal{N}(x;D,E)$ as stated in Theorem \ref{thm:thmfundamental}. Now write $DE=x^\psi$ where $\psi>0$.
Then $UV=x^{2/l-k\psi/l}$. We pick the integer $M$ such that
\[
	\frac{\log M}{\log x} = \frac{9}{8}(1+\delta)\psi(2/l-k\psi/l).
\]
Note that
\[
	\frac{9}{8}\psi(2/l-k\psi/l)\leq\frac{9}{8kl}<1.
\]
So, indeed with our choice of $M$ we have that $\log M\leq\log x$. This finishes the proof of Theorem \ref{thm:thmfundamental}.

%% file: PairsOfKFree.tex
\section{The Proof of Theorem \ref{thm:thm1}}
\subsection{Preliminaries} \label{sec:preliminaries}
First, we illustrate on how finding an asymptotic formula for $N_{k,h}(x)$ can be reduced to counting points on the algebraic variety $e^kv-d^ku=h$ inside a certain bounded box. In what follows all implied constants may depend on $h$ and $k$. First observe that by dyadic subdivision it is enough to show that
\[
	N_{k,h}(2x)-N_{k,h}(x)=c_{k,h}x+O_{\epsilon}(x^{\omega(k)+\epsilon}).
\]
Let
\[
	\xi(n)=\left(\prod_{p^k\divides n} p\right)
		   \left(\prod_{p^k\divides n+h} p\right).
\]
Then $\xi(n)=1$ if and only if $n$ and $n+h$ are $k$-free. Thus, we may deduce that
\begin{align*}
	N_{k,h}(2x)-N_{k,h}(x)	&= \sum_{x<n\leq 2x} \sum_{m\divides\xi(n)}\mu(m) \\
			&= \sum_{m=1}^\infty \mu(m) N(x;m),
\end{align*}
where
\[
	N(x;m)=\card{x<n\leq 2x: \xi(n)\equiv 0\pmod {m}}.
\]
Observe that the congruence $\xi(n)\equiv 0\pmod {p}$ has exactly $\rho_{k,h}(p)$ solutions modulo ${p^k}$. Thus, by an argument similar to the proof of the Chinese remainder theorem, the congruence
$\xi(n)\equiv 0\pmod {m}$ has exactly
\[
	\rho_{k,h}(m)=\prod_{p\divides m}\rho_{k,h}(p)
\]
solutions modulo $m^k$. Hence
\[
	N(x;m)=\rho_{k,h}(m)\left(\frac{x}{m^k}+O(1)\right).
\]
Note that for square-free $m$, we have that 
\[
	\rho_{k,h}(m)\leq 2^{\omega(m)}\ll m^\epsilon.
\]
Next, we introduce a parameter $y$ with $x^{1/k}\leq y\leq x^{2/(k+1)}$ which we shall pick in a moment. Now, we look at the small terms in the above sum corresponding to the values of $m$
with $m\leq y$. These terms contribute
\begin{align*}
	&\sum_{m\leq y} \mu(m)\rho_{k,h}(m)\left(\frac{x}{m^k}+O(1)\right)\\
	= &x\sum_{m\leq y} \frac{\mu(m)\rho_{k,h}(m)}{m^k}
		+ O\left(\sum_{m\leq y} \rho_{k,h}(m)\right) \\
	= &x\sum_{m=1}^\infty \frac{\mu(m)\rho_{k,h}(m)}{m^k}
		+O\left(x\sum_{m>y}\frac{\rho_{k,h}(m)}{m^k}
		+\sum_{m\leq y} \rho_{k,h}(m)\right)\\
	= &c_{k,h}x+O(x^{1+\epsilon}y^{1-k})+O(x^\epsilon y)\\
	= &c_{k,h}x+O(x^\epsilon y),
\end{align*}
where the last equality follows from $x^{1/k}\leq y$. To minimize the remaining error term, we pick $y=x^{1/k}$.
Thus, we can see that the values of $m$ with $m\leq x^{1/k}$ contribute our main term. Hence, we are left to consider the values of $m$ with $m>x^{1/k}$. For each such $m$ we write $m=de$ where $d^k\divides n$ and $e^k\divides n+h$. By dyadic subdivision, these values $d,e$ lie in $O((\log x)^2)$ boxes $D/2<d\leq D$, $E/2<e\leq E$
where $D,E\ll x^{1/k}$ and $DE\gg x^{1/k}$. Hence, for one such pair $D, E$ we must have
\begin{equation*}
	N_{k,h}(2x)-N_{k,h}(x)=c_{k,h}x+O(x^{1/k+\epsilon})
	+O(x^\epsilon\mathcal{N}(x;D,E)),
\end{equation*}
where $\mathcal{N}(x;D,E)$ is the number of elements in the set
\[
	\{ (d,e,u,v)\in\setN^4: D/2<d\leq D, E/2<e\leq E,
	x<d^ku+h=e^kv\leq 2x\}.
\]

\subsection{The proof of Theorem \ref{thm:thm1} for $k\geq 3$} \label{sec:proofkgeq3}
We apply Theorem \ref{thm:thmfundamental} with $l=1$, $U=x/D^k$ and $V=x/E^k$ to deduce that
\begin{equation}\label{bound1}
\mathcal{N}(x;D,E)\ll_{\epsilon}
	x^\epsilon\min\{(DEM)^{1/2}+D+E,(UVM)^{1/2}+U+V\}.
\end{equation}
We assume without loss of generality that $D\leq E$ and hence $V\leq U$. Note that $U/V\geq E/D\geq 1$. First we consider the case $M\geq U/V$. In this case, the estimate \eqref{bound1} becomes
\begin{equation}\label{bound2}
	\mathcal{N}(x;D,E)\ll_\epsilon x^\epsilon M^{1/2}\min(DE,UV)^{1/2}.
\end{equation}
We now set $DE=x^\psi$ so that $UV=x^{2-k\psi}$. Then
\[
	\frac{\log\mathcal{N}(x;D,E)}{\log x}\leq \epsilon+\frac{1}{2}\min(\psi,2-k\psi)+\frac{9}{16}\psi (2-k\psi)=:f_k(\psi),
\]
say. The function $f(\psi)$ takes its maximal value at $\psi=2/3$ if $k=2$ and at $\psi=13/(9k)$ if $k\geq 3$. Note that $f_2(2/3)=7/12$ and 
$f_k(13/(9k))=169/(144k)$ for $k\geq 3$. Note that
\[
	(DEM)^{1/2}+D+E\ll x^{\frac{169}{144k}} + x^{1/k}\ll x^{\frac{169}{144k}},
\]
no matter if $M\geq U/V$ does or does not hold. This completes the proof of Theorem \ref{thm:thm1} for $k\geq 3$. However, for $k=2$, we have
$169/(144k)\geq 7/12$ and a better argument is required. Thus, we shall employ the following trivial estimate, where we first sum over the pairs $e,u$ as follows:
\begin{align}
	\mathcal{N}(x;D,E) &\ll \mathop{\sum_{E\ll e\ll E}}_{U\ll u\ll U} 
	\#\{d: d\asymp D, e^2 v-d^2u= h\}\notag \\
	&\ll\mathop{\sum_{E\ll e\ll E}}_{U\ll u\ll U} 
	\#\{d: d\asymp D, d^2\equiv -hu^{-1}\pmod{e^2}\}\notag \\
	&\ll\mathop{\sum_{E\ll e\ll E}}_{U\ll u\ll U} 
	\left(\frac{D}{e^2}+1\right)
	\#\{d\mod{e^2}, d^2\equiv -hu^{-1}\pmod{e^2}\}\notag \\
	&\ll EU\left(\frac{D}{E^2}+1\right)x^\epsilon.\label{eq:eqtrivialeu1}
\end{align}
Note that 
\[
	EU\frac{D}{E^2}=\frac{x}{DE}\ll x^{1/2},
\]
Thus we may assume that $EU\geq x^{1/2}$ and, hence:
\[
	\frac{D^3}{E^3}=\frac{x^2}{(EU)^2ED}\ll\frac{x^2}{(x^{1/2})^3}=x^{1/2}.
\]
Thus, $\frac{D}{E}\ll x^{1/6}$ and $\frac{V}{U}\ll x^{1/3}$. By interchanging $E$ and $D$ and doing the same argument, we may deduce that
all of the quotients $\frac{D}{E}$, $\frac{E}{D}$, $\frac{V}{U}$, $\frac{U}{V}$ are $\ll x^{1/3}$. We may impose the condition $M\geq x^{1/3}\geq U/V$, so that the bound \eqref{bound2} does indeed hold. This shows that we can take the value $\omega(2)=7/12$ in Theorem \ref{thm:thm1}. But we can do better than this by considering points on lines contained in the three-fold $e^2v-d^2u=h$. The following section will illustrate this idea.

%% file: FurtherImprovements.tex
\subsection{Counting Points on Lines, Finishing the proof for $k=2$}\label{sec:FurtherImprovements}

In this section, we will conclude the proof of Theorem \ref{thm:thm1} by considering points on lines contained in the three-fold \eqref{eq:eq1}. For convenience, we will illustrate the proof when $h=1$ but with minor changes, the proof can be adapted to general $h$.
We begin by a similar procedure as in the proof of Theorem \ref{thm:thmfundamental}. We will pick $M\in [x^{1/2},x]$ and $A=B=1$ as in section \ref{sec:DetMethod}. Note that this will produce a $4\times 4$ matrix $\mathcal{M}$ where the $j$-th row is 
\[
	(\begin{array}{cccc} 
		1 & s_0(1+\alpha_j) & s_0^2(1+2\alpha_j+\beta_j) 
		  & s_0^3(1+\alpha_j)(1+2\alpha_j+\beta_j)
	\end{array}),
\]
with $\alpha_j\ll M^{-1}$ and $\beta_j\ll x^{-1}$. (Note that $\alpha_j^2\ll M^{-2}\leq x^{-1}$). Performing column operations, we get a matrix with $j$-th row 
\[
	s_0^6\cdot(\begin{array}{cccc} 
		1 & \alpha_j & \beta_j & \alpha_j (2\alpha_j+\beta_j)
	\end{array}).
\]
Recalling that $s_0\asymp D/E$, we may deduce that $\Delta\ll \frac{1}{M^3 x}\frac{D^6}{E^6}$. As before, we require $\abs{\Delta}\ll E^{-4}U^{-4}$
in order to deduce that $\Delta=0$. This is satisfied if $M>(xUV)^{1/3+\delta}$ where $\delta>0$ is arbitrarily small
and $x$ is large enough in terms of $\delta$. Setting as before $DE=x^\psi$ and $UV=x^{2-2\psi}$
we can see that $x^{1/2}\leq (xUV)^{1/3+\delta}\leq x$. This enables us to pick $M=(xUV)^{1/3+\delta}$ provided $\delta$ is small enough
and $x$ is large enough in terms of $\delta$.\bigskip\\
Hence, as before, the determinant method will produce an irreducible auxiliary polynomial $F(d,e;u,v)$. This polynomial
will be bilinear since we picked the monomials in $\mathcal{M}$ accordingly. That is, in section \ref{sec:Counting} the case that will occur is $a=b=1$.
Hence, \eqref{eq:eqG1} can be written as
\[
	d L_1(u,v)+ e L_2(u,v)=0,
\]
where $L_1(u,v)=c_1 u+c_2 v$ and $L_2(u,v)=c_3 u+c_4 v$ are linear forms with integral coefficients.
As in section \ref{sec:Counting}, we get $O(x^\delta)$ choices for an integer $q$ such that
\begin{equation}\label{eq:dqeq}
 e q=-L_1(u,v),\qquad d q=L_2(u,v).
\end{equation}
Plugging this into the equation $e^2v-d^2u=1$ gives us a Thue equation
\[
	G_3(u,v)=(L_1(u,v))^2v-(L_2(u,v))^2u=q^2,
\]
say. When $G_3$ is irreducible then this equation will only have a finite number of solutions (see Thue \cite{Thue}). Thus, the equation will produce $O(x^\delta)$ solutions $(u,v)$ for each $q$ except if $G_3$ is
splitting into three equal linear factors,
\[
	G_3(u,v)=\alpha (\alpha_1 v-\alpha_2 u)^3, 
\]
say where $\alpha,\alpha_1,\alpha_2\in\setZ$. Our aim is to show that the points under consideration corresponding to this case
actually lie on a line contained in the three-fold defined by \eqref{eq:eq1}.
Comparing coefficients we get the equations
\begin{align*}
	\alpha\alpha_2^3 &=c_3^2\\
	\alpha\alpha_1^3 &=c_2^2\\
	3\alpha\alpha_1^2\alpha_2 &=c_4^2-2c_1 c_2\\
	3\alpha\alpha_1\alpha_2^2 &=c_1^2-2c_3 c_4.
\end{align*}
The first two equations give $\alpha_1=(c_2^2/\alpha)^{1/3}$ and $\alpha_2=(c_3^2/\alpha)^{1/3}$ which turns the third and fourth
equation into
\begin{align}
	3c_2^{4/3}c_3^{2/3} &=c_4^2-2c_1 c_2,\label{eq:eqc1c2}\qquad\text{and}\\
	3c_2^{2/3}c_3^{4/3} &=c_1^2-2c_3 c_4 \label{eq:eqc3c4}
\end{align}
respectively. If $c_2=0$ then we may deduce from these equations that $c_4=0$ and hence $c_1=0$. Hence $L_1=0$ which gives a contradiction
since $qe\neq 0$. Thus, $c_2\neq 0$ and similarly $c_3\neq 0$. Using \eqref{eq:eqc1c2} and \eqref{eq:eqc3c4} we may deduce that $c_3^{2/3}(c_4^2-2c_1c_2)=c_2^{2/3}(c_1^2-2c_3c_4)$ which
may be written as 
\[
	(c_3^{1/3}c_4+c_3^{2/3}c_2^{2/3})^2=(c_2^{1/3}c_1+c_2^{2/3}c_3^{2/3})^2.
\]
After taking the square-root, this either lets us deduce directly that
$c_4^3 c_3=c_1^3 c_2$ or that 
\[
	c_3^{1/3}c_4+c_2^{1/3}c_1+2c_2^{2/3}c_3^{2/3}=0.
\]
Multiply this equation by $c_3^{2/3}$ and plug in the expression for $c_3 c_4$ given by \eqref{eq:eqc3c4} to get
$(c_1+c_2^{1/3}c_3^{2/3})^2=0$. Similarly, we may deduce from \eqref{eq:eqc1c2} that $(c_4+c_3^{1/3}c_2^{2/3})^2=0$.
This implies that in either case we have $c_4^3 c_3=c_1^3 c_2$. Plugging this value of $c_3$ into \eqref{eq:eqc1c2}
we get that
\[
	c_1=\frac{1}{\kappa}\frac{c_4^2}{c_2},\quad c_3=\frac{1}{\kappa^3}\frac{c_4^3}{c_2^2},
\]
where $\kappa\in\lbrace -1, 3\rbrace$. Consider the equation $c_4^3 c_3=c_1^3 c_2$. We set
$c_1=h\alpha$ and $c_4=h\beta$ where $h=(c_1,c_4)$. From this, we can deduce that $c_2=k\beta^3$ and $c_3=k\alpha^3$
for some integer $k$. Observe that $(h,k)=1$ since $(c_1,c_2,c_3,c_4)=1$ as $F$ is irreducible.
Considering the equation $c_4^2=\kappa c_1 c_2$ we see that $h=\kappa k\alpha\beta$ which implies $k=1$ and
\[
	c_1=\kappa\alpha^2\beta,\ c_2=\beta^3,\ c_3=\alpha^3,\ c_4=\kappa\alpha\beta^2.
\]
Our original Thue equation has now turned into
\[
	v=\frac{\alpha^2}{\beta^2} u + \frac{Q^2}{\beta^2},
\]
where $q=Q^3$ for $Q\in\setZ$ with $Q>0$. Substituting this expression for $v$ into \eqref{eq:dqeq}, we obtain
that $(d,e,u,v)$ must lie on the line
\[
	\ell:(d,e,u,v)=\left(\frac{\kappa\alpha}{Q},-\frac{\beta}{Q},0,\frac{Q^2}{\beta^2}\right)
				+u\left((1+\kappa)\frac{\alpha^3}{Q^3},-(\kappa+1)\frac{\alpha^2\beta}{Q^3},1,\frac{\alpha^2}{\beta^2}\right).
\]
Assume that the line does indeed have an integral point counted by $\mathcal{N}(x
;D,E)$. Pick $u_1$ to be the smallest
integer such that $(d_1,e_1,u_1,v_1)$ is on $\ell$ and counted by $\mathcal{N}(x;D,E)$. This gives us the equations
\begin{align}
	Q^2=v_1\beta^2-u_1\alpha^2\label{eq:eql1}\\
	\alpha\kappa Q^2+u_1 (1+\kappa)\alpha^3=d_1 Q^3\label{eq:eql2}\\
	-\beta Q^2-u_1(1+\kappa)\alpha^2\beta=e_1 Q^3.\label{eq:eql3}
\end{align}
We now consider the lines with $\kappa=-1$. In this case, we may deduce from \eqref{eq:eql2} and \eqref{eq:eql3}
that $Q=1$, $d_1=-\alpha$ and $e_1=-\beta$. Thus our points under consideration must lie on the line
\[
	(d,e,u,v)=(d_1,e_1,u_1,v_1)+\lambda (0,0,\beta^2,\alpha^2),
\]
where $\lambda\in\setZ$. Note that there are no parallel lines with $\kappa=-1$ having the same $\alpha$ and $\beta$ since
\eqref{eq:eql1} determines the pair $(u_1,v_1)$ modulo $(\beta^2,\alpha^2)$. The number of distinct lines with $\kappa=-1$
is therefore determined by the number of choices for the direction vectors, that is, by the number of choices for $\alpha$ and $\beta$.
We must have $\alpha\beta\ll (UV)^{1/2}\ll x^{1/2}$ which shows that the number of distinct lines with $\kappa=-1$ is $O(x^{1/2+\delta})$.
The number of points on each line is $\ll 1+ U/\beta^2\ll 1$ since $E\asymp e_1=-\beta$.
Thus, the total number of points counted be $\mathcal{N}(x;D,E)$ lying on lines with $\kappa=-1$ is $O(x^{1/2+\delta})$.\bigskip\\
Next, we consider the lines with $\kappa=3$. From \eqref{eq:eql2} and \eqref{eq:eql3} we may deduce that $(\beta d_1+\alpha e_1)Q=2\alpha\beta$.
Let $a=(Q,\alpha)$ and $b=(Q,\beta)$ so that $\alpha=a A, \beta=b B$ and $Q=a b Q_1$ and $(Q_1, A B)=1$ say.
After canceling $ab$ from the last equation, we may deduce that $A\divides d_1$ and $B\divides e_1$.
Set $d_1=A D_1$ and $e_1=B E_1$. Finally, \eqref{eq:eql1} implies that $a^2\divides v_1$ and $b^2\divides u_1$. Let
$v_1=a^2 V_1$ and $u_1=b^2 U_1$. Furthermore, we can deduce from \eqref{eq:eql1} that $(U_1,Q_1)=1$ and hence 
\eqref{eq:eql3} shows that $Q_1^2\divides 4 U_1 A^2$ and therefore $Q_1^2\divides 4$.
The direction vector of the line $\ell$ is now
\[
	\left(\frac{A^3}{b^3} C, -\frac{A^2 B}{a b^2} C, 1, \frac{a^2A^2}{b^2B^2}\right),
\]
where $C=4/Q_1^3\in\lbrace 4, 1/2\rbrace$. Note that $a,A,b,B$ are pairwise coprime since 
$d_1^2u_1$ and $e_1^2 v_1$ as well as $\alpha$ and $\beta$ are coprime. Thus, all lines lying on the surface \eqref{eq:eq1} can be written as 
\[
	(A D_1, B E_1, b^2U_1, a^2 V_1) +\lambda (a A^3 B^2 C, -b A^2B^3 C, ab^3 B^2, a^3 b A^2),
\]
where $\lambda$ goes through $\setZ$, $\setZ/2$ or $2\setZ$ and subject to the conditions
\begin{align}
	Q_1^2=V_1 B^2 - U_1 A^2\label{eq:eql4}\\
	b D_1 + a E_1 = 2 Q_1^{-1}\label{eq:eql5}\\
	U_1 A^2 Q_1 C - D_1 b Q_1 = -3\label{eq:eql6}.
\end{align}
Fix $Q_1$. Given a line $\ell$, the direction vector and hence $A,B,a,b$ are uniquely determined up to sign. Thus we also fix $A$ and $B$ and consider the number of points on lines with our fixed values of $Q_1, A$ and $B$.\bigskip\\
If $D_1',E_1',U_1',V_1'$ is one solution of \eqref{eq:eql4}, \eqref{eq:eql5}, \eqref{eq:eql6} then the other 
solutions $D_1,E_1,U_1,V_1$ are given by 
\begin{align*}
	D_1=D_1' +\mu a A^2 B^2\\
	E_1=E_1' -\mu  b A^2 B^2\\
	U_1=U_1' +\mu ab B^2/C\\
	V_1=V_1' +\mu ab A^2/C.
\end{align*}
Thus, the number of distinct lines with $\kappa=3$ and given $A,B$ is 
\[ 
\begin{split}
	&\ll\card{(a,b): a\ll D, b\ll E, ab^3\ll U, a^3b\ll V}\\
    &\ll\card{(a,b): ab\ll (UV)^{1/4}}\ll (UV)^{1/4+\delta}.
\end{split}
\]
The number of points on each line is
\[
	\ll\min\left\lbrace\frac{D}{ a A^3 B^2},\frac{E}{ b A^2 B^3}\right\rbrace\ll\frac{\min(D,E)}{A^2 B^2}.
\]
Thus, the total number of points on lines with $\kappa=3$ counted by
$\mathcal{N}(x;D,E)$ is
\[
	\sum_{A,B} (UV)^{1/4+\delta}\frac{\min(D,E)}{A^2 B^2}\ll\min(D,E)(UV)^{1/4+\delta}\ll x^{1/2+\delta}.
\]
Let 
\[
	\mathcal{N}_0(x;D,E)=\lbrace (d,e,u,v)\in \mathcal{N}(D,E): (d,e,u,v) \text{ is not on a line contained in \eqref{eq:eq1}}\rbrace.
\]
We have just shown that
\[
	\mathcal{N}_0(x;D,E)\ll_\delta x^\delta M\ll_\delta x^{2\delta+1-2\psi/3}.
\]
In section \ref{sec:proofkgeq3}, we have deduced the bound \eqref{bound2}. Hence,
\[
	\mathcal{N}_0(x;D,E)\leq \mathcal{N}(x;D,E)\ll_{\delta} x^{\delta+\min(\psi,2-2\psi)/2+\max(1/6,9\psi(1-\psi)/8)},
\]
where essentially $1/2\leq\psi\leq 1$. One can check that the worst value for the exponent occurs if $\psi\leq 2/3$
and then that the worst value for $\psi$ must satisfy 
\[
	1-\frac{2}{3}\psi=\frac{9}{8}\psi(1-\psi)+\frac{\psi}{2}.
\]
Hence the critical value for $\psi$ is $(55-\sqrt{433})/54= 0.6331\ldots$ and we may deduce that
$\mathcal{N}_0(x;D,E)\ll_\delta x^{3\delta+\omega}$  where $\omega=(26+\sqrt{433})/81\leq 0.5779\leq 7/12=0.5833\ldots$. By the above argument,
the points counted by $\mathcal{N}(x;D,E)$ but not being elements of $\mathcal{N}_0(x;D,E)$ contribute $O(x^{1/2+\delta})$. This finishes the proof of Theorem \ref{thm:thm1}.

%% file: TouplesOfKFree.tex
\section{The Proof of Theorem \ref{thm:thm2}}\label{sec:Corollary}
We are now turning to the proof of Theorem \ref{thm:thm2}. We define
\[
	\xi(n)=\prod_{i=1}^r\prod_{p^k\divides l_i(n)}p,
\]
and for $w>1$ let
\[
	\mathcal{P}(w)=\prod_{p<w}p
\]
be the product of primes $p<w$. Furthermore, for $z>1$ define
\[
	\mathcal{S}(z)=\mathop{\sum_{x<n\leq 2x}}_{(\xi(n),\mathcal{P}(z))=1} 1.
\]
That is, $\mathcal{S}(z)$ is the number of integers $n$ in the interval $(n,2n]$ so that $l_1(n),\ldots,l_r(n)$ all do not have any $k$-th power prime divisor $p^k$ with $p<z$. In particular, $N(2x)-N(x)=\mathcal{S}(O(x^{1/k}))$. We also define
\[
	\mathcal{S}_d(w) = \sum_{\substack{x<n\leq 2x \\ \xi(n)\equiv 0\pmod{d}\\ (\xi(n),\mathcal{P}(w))=1}} 1.
\]
Next, we will employ the following identity which is essentially Buchstab's identity. That is, for $1<w<z$ observe that
\begin{equation*}
	\mathcal{S}(z)=\mathcal{S}(w)-\sum_{w\leq p<z}\mathcal{S}_p(p).
\end{equation*}
Applying the identity twice, we obtain for $w>1$:
\[
	N(2x)-N(x)=\mathcal{S}(O(x^{1/k}))=\mathcal{S}(w)-\sum_{w\leq p\ll x^{1/k}}\mathcal{S}_p(w)+\sum_{w\leq q<p\ll x^{1/k}}\mathcal{S}_{pq}(q).
\]
First, we will estimate the sum $\sum_{p}\mathcal{S}_p(w)$. We split
the sum over $p$ in two parts, the first range will be $w\leq p< y$ and for the second range $y\leq p\ll x^{1/k}$ a trivial estimate will suffice. Similarly to section \ref{sec:preliminaries} we can deduce that:
\begin{equation}\label{eq:eqtrivialNd}
	N(x;d)=\card{x<n\leq 2x:\xi(n)\equiv 0\pmod d}=\rho(d)\left(\frac{x}{d^k}+O(1)\right).
\end{equation}
Thus,
\begin{align*}
	\sum_{y\leq p\ll x^{1/k}}\mathcal{S}_p(w) &\ll 
	\sum_{y\leq p\ll x^{1/k}}N(x;p)\ll x^{\epsilon}\sum_{y\leq p\ll x^{1/k}}
	\left(\frac{x}{p^k}+1\right)\\
	&\ll x^\epsilon (xy^{1-k}+x^{1/k})\ll x^{1/k+\epsilon}.
\end{align*}
To estimate the the sum over the remaining range  $w\leq p<y$ we will apply the following Fundamental Sieve Lemma due to Heath-Brown \cite{RHBprimesinshortintervals} which adapted to our purpose states as follows:
\begin{lemma}\label{lemm:RHBsieve}
	For $z>1$ and $w>1$ we have that
	\[
		\sum_{d\divides (\xi(n),\mathcal{P}(w))}\mu(d)
		=\sum_{\substack{d\divides (\xi(n),\mathcal{P}(w))\\d<z}}\mu(d)
		+O\left(\sum_{\substack{
		d\divides (\xi(n),\mathcal{P}(w))\\
		z\leq d<zw}}1\right).
	\]
\end{lemma}
Observe that a prime $p>w>1$ and some $d\divides\mathcal{P}(w)$ both divide $\xi(n)$
if and only if $pd$ divides $\xi(n)$.
Thus, we pick some parameter $z_p>1$ which we will determine later and we split $\mathcal{S}_p(w)$ using the Lemma as follows:
\[
\begin{split}
	\mathcal{S}_p(w) &= \sum_{\substack{x<n\leq 2x \\ \xi(n)\equiv 0\pmod{p}}} 
	\sum_{d\divides (\xi(n),P(w))}\mu(d) \\	
	&=\sum_{\substack{x<n\leq 2x \\ \xi(n)\equiv 0\pmod{p}}} \sum_{
	\substack{d\divides (\xi(n),\mathcal{P}(w))\\d<z_p}}\mu(d)
	+O\left(\sum_{\substack{x<n\leq 2x \\ \xi(n)\equiv 0\pmod{p}}}\sum_{\substack{
		d\divides (\xi(n),\mathcal{P}(w))\\
		z_p\leq d<z_pw}}1\right)\\
	&=\sum_{\substack{d\divides\mathcal{P}(w)\\ d<z_p}}\mu(d)N(x;pd)+
	O\left(\sum_{\substack{d\divides\mathcal{P}(w)\\ 
	z_p\leq d<z_pw}}N(x;pd)\right)\\
	&=S_1(p)+O(S_2(p)),
\end{split}
\]
say. For $S_2(p)$ we obtain
\[
	S_2(p) \ll x^\epsilon \sum_{z_p\leq d<z_p w} \left(\frac{x}{d^kp^k}+1\right)
	\ll x^\epsilon\left(\frac{x}{p^k z_p^{k-1}}+z_p w\right).
\]
We will pick $z_p=p^{-1} (x/w)^{1/k}$ to minimize this error term.
Here we set $y=(x/w)^{1/k}$ to ensure that $z_p>1$. Note that $y=pz_p$. Using the fact that
$\sum_{p\leq x} p^{-k}\ll x^{\epsilon}$ we can deduce that
\[
	\sum_{w\leq p< y}S_2(p)\ll x^{1/k+\epsilon}w^{1-1/k}.
\]
Next, we use the trivial estimate \eqref{eq:eqtrivialNd} again to conclude that
\begin{align*}
	-\sum_{w\leq p< y}S_1(p) &=x\sum_{w\leq p< y}\sum_{
	\substack{d\divides\mathcal{P}(w) \\ pd< y}}
	\frac{\mu(pd)\rho(pd)}{(pd)^k}+O(yx^{\epsilon}).
\end{align*}
The double sum equals
\begin{align}
	x\sideset{}{'}\sum_{d<y}\frac{\mu(d)\rho(d)}{d^k}
	&=x\sideset{}{'}\sum_{d=1}^\infty\frac{\mu(d)\rho(d)}{d^k}
	+O(x^{1+\epsilon}y^{1-k}),\label{eq:eqc1}
\end{align}
where the $\sideset{}{'}\sum$ restricts the sum to those $d\divides \mathcal{P}(x^{1/k})$ with exactly one exceptional prime divisor $p\divides d$ such that $p> w$. Thus, we have shown that
\[
	-\sum_{w\leq p\ll x^{1/k}}\mathcal{S}_p(w)=c^{(1)}x+O(x^{1/k+\epsilon}w^{1-1/k}),
\]
where $c^{(1)}$ is the constant from \eqref{eq:eqc1}, the sum over those $d$ having exactly one large prime divisor $p>w$. Next, we consider the sum $\mathcal{S}(w)$. Similarly to the above argument, we may apply Lemma \ref{lemm:RHBsieve} for some $z>1$ to obtain
\begin{equation}\label{eq:eqSw}
	\mathcal{S}(w)= 
	\sum_{\substack{d\divides\mathcal{P}(w)\\d<z}} \mu(d)N(x;d)+
	O\left(\sum_{\substack{ d\divides P(w)\\z\leq d<zw}} N(x;d)\right).
\end{equation}
Again, we use a trivial estimate for the second sum which yields an error term \\$O(x^{1/k+\epsilon}w^{1-1/k})$ provided we chose $z=(x/w)^{1/k}$ optimally. The first sum in \eqref{eq:eqSw}
is
\[
	x\sum_{\substack{d\divides\mathcal{P}(w)\\d<z}}\frac{\mu(d)\rho(d)}{d^k}+O(z).
\]
Thus, we have shown that
\[
	\mathcal{S}(w)=c^{(0)}x+O(x^{1/k+\epsilon}w^{1-1/k}),
\]
where
\[
	c^{(0)}=\sum_{d\divides \mathcal{P}(w)}\frac{\mu(d)\rho(d)}{d^k}
\]
is the sum over those $d$ having no large prime divisors $p>w$. Recall that the overall main
term is $cx=x\sum_{d=1}^\infty \mu(d)\rho(d)d^{-k}$. The sum
$c-c^{(0)}-c^{(1)}$ is the sum over those $d$ having at least 2 distinct prime divisors $>w$ and thus we have
\[
	c^{(0)}+c^{(1)}=c+O\left(x^\epsilon\sum_{d>w^2}\frac{1}{d^k}\right)=
	c+O(x^{\epsilon}w^{2-2k}).
\]
We minimize the error terms by choosing $w=x^{1/(2k+1)}$ so that
both our error terms $O(x^{1+\epsilon}w^{2-2k})$ and
$O(x^{1/k+\epsilon}w^{1-1/k})$ become $O(x^{3/(2k+1)+\epsilon})$. Thus, we have
shown
\[
	\mathcal{S}(w)-\sum_{w\leq p\ll x^{1/k}}\mathcal{S}_p(w) =cx+
	O(x^{3/(2k+1)+\epsilon}),
\]
and it remains to find a bound for the sum $\sum_{w\leq q<p\ll x^{1/k}}\mathcal{S}_{pq}(q)$. First we consider the terms corresponding to those prime pairs $q<p$ with $pq\ll x^{1/k}$. A trivial estimate suffices to yield the bound
\[
	\sum_{\substack{w\leq q<p\ll x^{1/k}\\ pq\ll
	x^{1/k}}}\mathcal{S}_{pq}(q)
	\ll \sum_{\substack{w\leq q<p\ll x^{1/k}\\ pq\ll x^{1/k}}}N(x;pq)
	\ll x^\epsilon\sum_{\substack{w\leq q<p\ll x^{1/k}\\ pq\ll
	x^{1/k}}}\left(
	\frac{x}{p^kq^k}+1\right)\ll x^{3/(2k+1)+\epsilon}.
\]
For the values with $pq\gg x^{1/k}$ observe that
\[
	\mathcal{S}_{pq}(q)\ll N(x;pq)=\#\{x<n\leq 2x:
	p^k\divides l_i(n), q^k\divides l_j(n)\text{ for some } i\neq j\}.
\]
(The case $i=j$ cannot occur since $pq\gg x^{1/k}$  for any suitable implied constant). Thus, we can fix a particular $i$ and $j$ and conclude that
\[
\begin{split}
	\sum_{\substack{w\leq q<p\ll x^{1/k}\\ pq\gg x^{1/k}}}
	\mathcal{S}_{pq}(q)
	\ll_r \#\mathcal{K},
\end{split}
\]
where
\[
	\mathcal{K}=\{(p,q,u,v): w\leq q<p\ll x^{1/k},pq\gg x^{1/k},
	p^ku=l_i(n), q^kv=l_j(n)\}.
\]
Without loss of generality we may write $l_i(n)=a_1 n+b_1$ and
$l_j(n)=a_2 n+b_2$ so that we are left with the Diophantine equation
$a_1 q^kv-a_2 p^ku=a_1 b_2-a_2 b_1=h\neq 0$ say. Note that
$(a_1 q^k v, a_2 p^k u)=(a_1 v,a_2 u)$ since neither $p$ nor $q$ divide $h$ since $h$ is $O(1)$ and $p$ and $q$ are $\gg x^{1/(2k+1)}$.
Each common divisor of $a_1v$ and $a_2u$ is a divisor of $h$ and thus is $O(1)$.
Since also $a_1$ and $a_2$ are $O(1)$ we can reduce our problem to $d(h)=O(1)$ equations of the form \eqref{eq:eq1} with the additional constraint that $(du,ev)=1$ and $de\gg x^{1/k}$. Thus,
we have reduced the problem to the case we dealt with in Theorem $\ref{thm:thm1}$. Thus, this gives an error term $O(x^{\max\{\omega(k),3/(2k+1)\}+\epsilon})$ for the asymptotic formula in Theorem \ref{thm:thm2}. Note that $\omega(k)\leq 3/(2k+1)$ for all $k\geq 2$ which concludes the proof.

%% file: ConsecutiveSqFull.tex
\section{The Proof of Theorem \ref{thm:thm3}}
We are now turning to the problem of consecutive square-full integers. Recall that here, $N(x)$ is the number of integers $n\leq x$ such that both $n$ and $n+1$ are square-full. Observe that every square-full integer $n$ can uniquely be written as
$n=a^2b^3$ with $\mu^2(b)=1$. Thus, we have
\[
	N(2x)-N(x)\ll\card{(d,e,u,v)\in\setN^4: x<d^3u^2=e^3v^2-1\leq 2x}.
\]
As in the proof of Theorem \ref{thm:thm1} we can now split the ranges of $d$ and $e$ into $O(x^\epsilon)$ boxes with $D/2<d<D$ and $E/2<e\leq E$ where $D,E\ll x^{1/3}$. For one such box we then get
\[
	N(2x)-N(x)\ll x^\epsilon\mathcal{N}(x;D,E),
\]
where 
\[
	\mathcal{N}(x;D,E)=\card{(d,e,u,v)\in\setN^4: x<d^3u^2=e^3v^2-1\leq 2x,
	d\asymp D, e\asymp E}.
\]
Thus, in order to prove Theorem \ref{thm:thm3}, it remains to show that
\[
	\mathcal{N}(x;D,E)\ll_\epsilon x^{29/100+\epsilon}.
\]
For convenience, we will set $y=x^{29/100}$.
We may now apply Theorem \ref{thm:thmfundamental} with $k=3, l=2$ and $U:=x^{1/2}D^{-3/2}$, $V:=x^{1/2}E^{-3/2}$ so that $v\asymp V$ and $u\asymp U$.
We get the bound 
\begin{equation}\label{eq:eqxxx}
	\mathcal{N}(x;D,E)\ll_{\epsilon}
	x^\epsilon\min\{(DEM)^{1/2}+D+E,(UVM)^{1/2}+U+V\},
\end{equation}
where the value of $M$ is as stated in Theorem \ref{thm:thmfundamental}. We will further impose the condition $M\geq x^{9/50}$. We will see below that this will not make our bound worse. The situation we are now in is different compared to the situation in section \ref{sec:proofkgeq3}, where we used the estimate $D+E\ll x^{1/k}$. In our situation here this would produce an estimate $\mathcal{N}(x;D,E)\ll_\epsilon x^{1/3+\epsilon}$ which is worse than our anticipated bound. And without further restricting the ranges of $D$ and $E$, it will not be possible to improve this bound. The results by Estermann and Thue allow us to assume further that $\min(DE,UV)\geq y$. This is however still not enough to improve upon the estimate $\mathcal{N}(x;D,E)\ll_\epsilon x^{1/3+\epsilon}$. Thus, we shall employ an estimate similar to \eqref{eq:eqtrivialeu1}, namely:
\begin{align}
	\mathcal{N}(x;D,E) &\ll_\epsilon 
	EU\left(\frac{D}{E^3}+1\right)x^\epsilon\label{eq:eqtrivialeu2}.
\end{align}
We will now distinguish two cases depending on whether $DE^{-3}\ll 1$ does or does not hold. If $DE^{-3}\ll 1$ for some suitable constant, then we may assume $EU\geq y$ by the above estimate. Note that by definition of $U$,
\[
	y^5\leq (DE)(EU)^4=\left(\frac{E}{D}\right)^5 x^2,
\]
so that 
\[
	\frac{V}{U}=\left(\frac{D}{E}\right)^{3/2}\leq \frac{x^{3/5}}{y^{3/2}}
	=x^{33/200}< x^{9/50}\leq M.
\]
Next, we consider the case $DE^{-3}\gg 1$. Here, the estimate \eqref{eq:eqtrivialeu2} is too weak for our purposes and we need a stronger approach.
The idea is to fix a pair $(e,u)$ and estimate the number of solutions $(d,e,u,v)$ counted by $\mathcal{N}(x;D,E)$. Let $d_1,\ldots,d_\mu$ be the solutions to the congruence $d^3\equiv -u^{-2} \pmod {e^3}$.
In particular, $\mu\ll e^\epsilon\ll x^\epsilon$ and we get the estimate
\[
	\mathcal{N}(x;D,E)\ll\sum_{e\sim E}\sum_{u\sim U}\sum_{i=1}^\mu 
	\#\{(a,v): v\sim V, a\ll D/E^3, e^3v^2-(d_i+ae^3)^3u^2=1\}.
\]
Thus, it suffices to estimate the number of points $(a,v)$ that satisfy
\begin{equation}\label{eq:eqpq}
	p(v)-q(a)-1=0,
\end{equation}
where $p$ is a polynomial of degree $2$ and $q$ is a polynomial of degree $3$.\smallskip\\
We will now prove that that the left-hand side of \eqref{eq:eqpq} is absolutely irreducible. It is enough to show that the polynomial $f(S,T)=S^2-T^3-1$ is absolutely irreducible. Assume that $f(S,T)=g(S,T)h(S,T)$ in some finite extension of $\setQ$. In particular, $X^6-Y^6-1=g(X^3,Y^2)h(X^3,Y^2)$. We may homogenize the equation to get
\begin{equation}\label{eq:eqx6y6z6}
	X^6-Y^6-Z^6=G(X,Y,Z)H(X,Y,Z)
\end{equation}
for some polynomials $G,H$. Now let $(X,Y,Z)$ be a nonzero point such that $G(X,Y,Z)=H(X,Y,Z)=0$. The gradient of the right-hand side of \eqref{eq:eqx6y6z6} vanishes whereas the gradient of the left-hand side is $(6X^5,-6Y^5,-6Z^5)$ which implies that $X=Y=Z=0$. Contradiction. Thus, the equation \eqref{eq:eqpq} must be absolutely irreducible.\smallskip\\
Next, we will proceed to apply Theorem 15 of Heath-Brown \cite{RHBratptsonalgvar} to get an upper bound for the number of pairs $(v,a)$ satisfying \eqref{eq:eqpq}. Using the notation of Heath-Brown's theorem we will set $n:=2$, $B_1=V$ and $B_2=DE^{-3}+1$ so that indeed $v\leq B_1$ and $a\leq B_2$ with $B_1,B_2\geq 1$. Note that 
$T=\max(B_1^2,B_2^3)\geq B_1^2$ so that the points $(v,a)$ satisfying \eqref{eq:eqpq} lie on at most $k\ll_\epsilon x^\epsilon B_2^{1/2}$ auxiliary curves. Thus, using B\'{e}zout's Theorem, we may deduce that the number of points $(v,a)$ under consideration is $\ll_\epsilon x^\epsilon B_2^{1/2}$.
Thus, we get the estimate
\[
	\mathcal{N}(D,E)\ll_\epsilon x^\epsilon EU\left(\frac{D}{E^3}+1\right)^{1/2}\ll\frac{x^{1/2+\epsilon}}{DE^{1/2}},
\]
where the last estimate is because $1\ll D/E^{3}$. Hence, we may assume that 
\[
	x^{1/2}D^{-1}E^{-1/2}\geq y.
\]
Thus,
\[
	y^7\leq\left(\frac{x^{1/2}}{DE^{1/2}}\right)^4(DE)^3=x^2\frac{E}{D},
\]
so that
\[
	\frac{D}{E}\leq\frac{x^2}{y^7}\leq 1.
\]
Thus, we have shown that in all cases, $V/U\leq M$. By interchanging the roles of $D$ and $E$, we may similarly prove that $U/V\leq M$ and hence we conclude that
\begin{equation*}\label{eq:eqmaxdeeduvvu}
	\max(D/E,E/D,U/V,V/U)\leq M.
\end{equation*}
By considering \eqref{eq:eqxxx} again, we may now deduce that
\[
	\mathcal{N}(x;D,E)\ll_{\epsilon}
	x^\epsilon M^{1/2}\min\{DE,UV\}^{1/2}.
\]
We will set $DE=x^\psi$ so that $UV=x^{1-3\psi/2}$. Observe that
\[
	\frac{\log\mathcal{N}(x;D,E)}{\log x}
	\leq \delta+\frac{1}{2}\min(\psi,1-3\psi/2)+\frac{1}{2}
	\max\left\lbrace\frac{9}{8}\psi (1-3\psi/2),\frac{9}{50}\right\rbrace:=f(\psi),
\]
say. The function $f(\psi)$ takes its maximum at $\psi=2/5$ and $f(2/5)=29/100$. This completes the proof of Theorem \ref{thm:thm3}.

%% file: PellEquation.tex
\section{The Proof of Theorem \ref{thm:thmPell1} and Corollary \ref{cor:corPell}}
In this proof, all implied constants may depend on $\alpha$.
Our initial argument is the same as in \cite{FouvryJouve}. We want to find $r=r(\alpha)$ as small
as possible so that the estimate $S(X,\alpha)\ll_\epsilon X^{r+\epsilon}$ holds. As illustrated in \cite{FouvryJouve}, we may assume that $\alpha\geq 1/2$. The estimate \eqref{eq:eqtrivFouvry}
allows us to assume that $r\geq 1/2$. As in \cite{FouvryJouve}, we obtain the estimate
$S(X,\alpha)\ll\mathcal{N}(X;D,E)$, where $\mathcal{N}(X;D,E)$ is the number of elements in the set
\[
	\#\{(d,e,u,v):d\sim D, e\sim E, u\sim U, v\sim V, e^2v-d^2u=h\},
\]
where $h\in\{\pm 1,\pm 2\}$ and $D,E,U,V,X$ are positive real numbers such that $E^2V\asymp D^2U$, $UV\asymp X$ and
$DE\ll X^\alpha$. In particular,
\begin{equation}\label{eq:eqforuv}
	V\asymp X^{1/2}D/E\quad\text{and}\quad U\asymp X^{1/2}E/D.
\end{equation}
We want to apply Theorem \ref{thm:thmfundamental} with $x=E^2V$, $k=2$ and $l=1$.
Note the trivial estimate
\[
	\mathcal{N}(X;D,E) \ll DE\left(\frac{V}{D^2}+1\right)\ll X^{1/2}+DE.
\]
Hence, if $DE\ll X^r$ then this estimate suffices. Thus, we may assume that $DE\geq X^{r}$. In particular, this shows that $DE\geq X^{1/2}$ and from $x^2\asymp (DE)^2X$ we can then conclude that 
$DE\gg x^{1/2}$. So, Theorem \ref{thm:thmfundamental} does indeed apply. We get the estimate
\[
	\mathcal{N}(X;D,E)\ll_{\epsilon}
	X^\epsilon \min\{(XM)^{1/2}+U+V,(DEM)^{1/2}+D+E\},
\]
where $D$, $E$, $U$ and $V$ are as in \eqref{eq:eqforuv} and
\[
	\log M =\frac{9}{8}\frac{\log(DE)\log X}{\log(DEX^{1/2})},
\]
since $x=E^2V\asymp DEX^{1/2}$ and $UV\asymp X$. As in the proof of Theorem \ref{thm:thm1}, we get the estimate
\[
	\mathcal{N}(X;D,E) \ll_\epsilon X^\epsilon EU\left(\frac{D}{E^2}+1\right)
	=X^\epsilon(X^{1/2}+EU).
\]
Note that $EU=X^{1/2}E^2/D$. Hence we may assume that
$E^2/D\geq X^{r-1/2}$ and similarly we assume that $D^2/E\geq X^{r-1/2}$. Without loss of generality, let $D\leq E$ so that $V\leq U$. First, we consider the case when $M<U/V=(E/D)^2$. In this case,
\[
	(XM)^{1/2}+U+V<X^{1/2}E/D+U+V\ll U.
\]
From $DE\ll X^\alpha$ and $D^2/E\geq X^{r-1/2}$, it follows that
$(E/D)^{3/2}\ll X^{(\alpha+1)/2-r}$, and hence
\[
	\mathcal{N}(X;D,E) \ll_\epsilon X^\epsilon U
	=X^{1/2+\epsilon}\frac{E}{D}
	\ll X^{\frac{1}{3}\alpha+\frac{5}{6}-\frac{2}{3}r+\epsilon}\ll X^{r+\epsilon},
\]
provided $r>1/2+\alpha/5$. This concludes the case $M<U/V$. So, we may now consider the case $M\geq U/V$.
Note that $M\geq U/V=(E/D)^2\geq E/D$ which yields the estimate
\[
	\mathcal{N}(X;D,E)\ll_{\epsilon} X^\epsilon M^{1/2}\min\{X,DE\}^{1/2}.
\]
We will set $DE=X^\psi$ so that essentially $\psi\leq \alpha$. Then
\[
	\frac{\log M}{\log X}=\frac{9}{8}\frac{\psi}{\psi+1/2}=:f(\psi),
\]
say. The function $f(\psi)$ is increasing for $\psi>0$ and hence $f(\psi)\leq f(\alpha)$. 
This shows that
\[
	\frac{\log\mathcal{N}(X;D,E)}{\log X}\leq
	 O_\epsilon(1)+\frac{1}{2}f(\alpha)+\frac{1}{2}\min(1,\alpha). 
\] 
This concludes the case $M\geq U/V$. We have therefore proved Theorem \ref{thm:thmPell1}. Corollary \ref{cor:corPell} is an easy consequence of Theorem \ref{thm:thmPell1}.